\documentclass[12pt,a4paper]{article}
\usepackage[cp1251]{inputenc}
\usepackage[T2A]{fontenc}
\usepackage[english]{babel}

\usepackage{amsmath,amsfonts,amssymb,mathrsfs,amsopn,indentfirst,amscd}

\usepackage{graphicx}
\usepackage{amsthm}
\usepackage{hyperref}

\usepackage{amssymb}
\usepackage{amsmath}

\newcommand{\ver}{{\rm ver}}
\newcommand{\vo}{{\rm vol}}

\newcommand{\lin}{{\rm lin}}

\newcommand{\conv}{{\rm conv}}

\newtheorem*{corollary*}{Corollary}

\begin{document}

\title{On a Geometric Approach to the Estimation\\ of Interpolation Projectors}

\author{Mikhail Nevskii\footnote{ Department of Mathematics,  P.G.~Demidov Yaroslavl State University,\newline Sovetskaya str., 14, Yaroslavl, 150003, Russia, \newline
               mnevsk55@yandex.ru,  orcid.org/0000-0002-6392-7618 },
\quad Alexey Ukhalov\footnote{ Department of Mathematics,  P.G.~Demidov Yaroslavl State University,\newline Sovetskaya str., 14, Yaroslavl, 150003, Russia, \newline
               alex-ukhalov@yandex.ru,  orcid.org/0000-0001-6551-5118 }
               }
      
\date{July 25, 2023}
\maketitle

\begin{abstract}

\smallskip
Suppose
$\Omega$ is  a closed bounded subset of
${\mathbb R}^n,$ $S$ is an $n$-dimensional non-degenerate simplex,
$\xi(\Omega;S):=\min \left\{\sigma\geq 1: \, \Omega\subset \sigma S\right\}$. Here $\sigma S$ is the result of 
homothety of $S$  with respect to the center of gravity with coefficient $\sigma$.
Let
$d\geq n+1,$
$\varphi_1(x),\ldots,\varphi_d(x)$ be linearly independent
monomials in $n$ variables, 
$\varphi_1(x)\equiv 1,$ $\varphi_2(x)=x_1,\
\ldots, \ \varphi_{n+1}(x)=x_n.$  Put $\Pi:=\lin(\varphi_1,\ldots,\varphi_d).$
The~interpolation projector $P: C(\Omega)\to \Pi$
with a set of nodes $x^{(1)},\ldots, x^{(d)}$ $ \in \Omega$ is~defined
by equalities
$Pf\left(x^{(j)}\right)=f\left(x^{(j)}\right).$
Denote by
$\|P\|_{\Omega}$
the norm of~$P$ as an operator from $C(\Omega)$ to
$C(\Omega)$ .
Consider the mapping
$T:{\mathbb R}^n\to {\mathbb R}^{d-1}$ \linebreak of the form
$T(x):=(\varphi_2(x),\ldots,\varphi_d(x)).
$
We have the following  inequalities:
$
\frac{1}{2}\left(1+\frac{1}{d-1}\right)\left(\|P\|_{\Omega}-1\right)+1$ $
\leq \xi(T(\Omega);S)\leq \frac{d}{2}\left(\|P\|_{\Omega}-1\right)+1.
$
Here $S$ \linebreak  is~the~$(d-1)$-dimensional simplex with vertices $T\left(x^{(j)}\right).$
We discuss
this and~other relations for polynomial interpolation of functions continuous \linebreak  on~a~segment.
The results of numerical analysis are presented.

\smallskip 



\end{abstract}

\section{Main Definitions and Relations}\label{nev_s1}

Let  
$\Omega$ be a closed bounded subset of 
${\mathbb R}^n$, $n\in {\mathbb N}$.
By
$C(\Omega)$ we mean the space of~continuous functions 
 $f:\Omega\to {\mathbb R}$  with the uniform norm
$$\|f\|_{C(\Omega)}:=\max_{x\in \Omega} |f(x)|.$$
For a nondegenerate simplex $S\subset {\mathbb R}^n$, by   $\sigma S$  denote the homothetic copy of $S$ with center of homothety at the center of gravity of $S$ and ratio $\sigma$.
By definition, put 
$\xi(\Omega;S):=\min \left\{\sigma\geq 1: \, \Omega\subset \sigma S\right\}$.
The inclusion  $\Omega\subset S$ is equivalent to the equality
$\xi(\Omega;S)=1.$ We call $\xi(\Omega;S)$ {\it the absorption coefficient of the set $\Omega$ by
the simplex $S$.} Further $\ver(S)$ is a set of vertices of $S$.

Assume $ d\in{\mathbb N},$ 
$d\geq n+1;$ 
$\varphi_1(x),\ldots,\varphi_d(x)$ are linearly independent functions, namely  monomials in $n$ variables having the form 
$x^\alpha=x_1^{\alpha_1}\ldots  x_n^{\alpha_n}.$
Here $x=(x_1,\ldots,x_n)
\in {\mathbb R}^n$,  $\alpha=(\alpha_1,\ldots,\alpha_n)$ $\in {\mathbb Z}_+^{n}$.
We suppose that
$\varphi_1(x)\equiv 1,$  $\varphi_2(x)=x_1, \ 
\ldots, \ \varphi_{n+1}(x)=x_n.$
By the $d$-dimensional space of polynomials in $n$ variables we mean the set
$\Pi:=\lin(\varphi_1,\ldots,\varphi_d).$
Note the important cases: 
$\Pi=\Pi_k\left({\mathbb R}^n\right)$ --- 
the space of polynomials of general degree $\leq k$\, $(k\in {\mathbb N})$  and
$\Pi=\Pi_\alpha\left({\mathbb R}^n\right)$ --- 
the~space of polynomials of degree $\leq \alpha_i$
in~$x_i$ \  $(\alpha\in{\mathbb N}^n).$

 A collection of points $x^{(1)},\ldots, x^{(d)} \in \Omega$ is called
{\it an admissible set of nodes } for~interpolation of functions from $C(\Omega)$ with the use of polynomials 
from $\Pi$  iff $\Delta:=\det({\bf A})\ne 0.$ Here
$\bf A$ is $(d\times d)$-matrix
 $${\bf A} :=
\left( \begin{array}{cccc}
1&\varphi_2\left(x^{(1)}\right)&
\ldots&\varphi_d\left(x^{(1)}\right)\\

\vdots&\vdots&\vdots&\vdots\\ 

1&\varphi_2\left(x^{(d)}\right)&
\ldots&\varphi_d\left(x^{(d)}\right)
\end{array}
\right).
$$  
The interpolation projector  $P: C(\Omega)\to \Pi$ related to this set of nodes is defined by~equalities
$Pf\left(x^{(j)}\right)=
%f_j:=
f\left(x^{(j)}\right), \ j=1,\ldots, d.$
An analogue of the Lagrange interpolation formula is the representation
\begin{equation}\label{Lagr_formula}
Pf(x)=\sum_{j=1}^d f\left(x^{(j)}\right)\lambda_j(x),
\quad
\lambda_j(x):=\frac{\Delta_j(x)}{\Delta}, 
\end{equation}
where $\Delta_j(x)$ is the determinant which appears from $\Delta$ by replacing the $j$th
row \linebreak with the row 
$\bigl(\varphi_1(x),\ldots,\varphi_d(x)\bigr).$ 
Polynomials $\lambda_j\in \Pi$ have the property \linebreak 
$\lambda_j\left(x^{(k)}\right)$ $=$ 
$\delta_j^k.$ Their coefficients with respect  to  the basis  $\varphi_1,$ $\ldots,$ $\varphi_d$ 
form the~columns of  ${\bf A}^{-1}.$
We call
$\lambda_j$ {\it the basic Lagrange polynomials corresponding to the  projector $P$}.

 In the case $\Pi=
\Pi_1\left({\bf R}^n\right),$
$d=n+1$, and
$\varphi_j(x)=x_{j-1}$ $(j=2,\ldots,d)$ polynomials
$\lambda_j$ are also called {\it the basic Lagrange polynomials of the simplex $S$} with vertices at~the~interpolation nodes.
If $\Omega\not\subset S$, then
\begin{equation}\label{xi_Lagr_pol_formula}
\xi(\Omega;S)=(n+1)  \max_{1\leq k\leq n+1} \max_{x\in \Omega} (-\lambda_k(x))+1.
\end{equation}
For convex $\Omega$, equality \eqref{xi_Lagr_pol_formula} is proved in  \cite{nevskii_monograph}; 
the proof in general case is just the same. 
Notice that  $\xi(\conv(\Omega);S)=\xi(\Omega;S).$

Further we consider only admissible sets of nodes and those $\Omega$ which contain such a set.

Denote by
$\|P\|_{\Omega}$ 
the norm of $P$ as an operator from $C(\Omega)$ in
$C(\Omega)$ .
From \eqref{Lagr_formula}, it~follows that
\begin{equation}\label{norm_P_by_basic_polynomials}
\|P\|_{\Omega}=  
 \max\limits_{x\in \Omega}
\sum_{j=1}^d |\lambda_j(x)|.
\end{equation}

The equalities $\lambda_j(x)=\dfrac{\Delta_j(x)}{\Delta}$ are equivalent to the matrix relation
 \begin{equation}\label{matrix_eq_for_bary_coords}
 \left ( \begin{array}{cccc}
1&1&\ldots&1\\
\varphi_2\left(x^{(1)}\right)
&\varphi_2\left(x^{(2)}\right)&\ldots&\varphi_2\left(x^{(d)}\right)\\
\vdots&\vdots& &\vdots\\
\varphi_d\left(x^{(1)}\right)
&\varphi_d\left(x^{(2)}\right)&\ldots&\varphi_d\left(x^{(d)}\right)\\
\end{array}
\right )
\left ( \begin{array}{c} \lambda_1\\\vdots\\ \lambda_{d}
\end{array} \right )
= \left ( \begin{array}{c} 1\\\varphi_2(x)\\\vdots\\\varphi_d(x)
\end{array} \right ).
\end{equation}
Let us introduce the mapping 
$T:{\mathbb R}^n\to {\mathbb R}^{d-1}$ defined in the way
$$y=T(x):=(\varphi_2(x),\ldots,\varphi_d(x))=
(x_1,\ldots,x_n,\varphi_{n+1}(x),\ldots,\varphi_d(x)).
$$
We will consider $T$ on the set $\Omega.$ 
The choice of the first monomials  $\varphi_j(x)$ implies the~invertibility of  $T.$
%Как обычно, $T(\Omega)$ обозначает образ $\Omega$ пpи отобpажении $T.$
Denote $y^{(j)}:=T\left(x^{(j)}\right)$. Relation  \eqref{matrix_eq_for_bary_coords} gives
$$
\lambda_1y^{(1)}+\ldots+\lambda_dy^{(d)}=y, \quad \sum_{j=1}^{d}\lambda_j = 1.
%\eqno{(1.5)}
$$
This means that the numbers $\lambda_j(x)$ are the barycentric coordinates of the point
$y=T(x)$ with respect to the $(d-1)$-dimensional simplex with vertices $y^{(j)}$.
Hence,
\begin{equation}\label{norm_P_by_bary_coordinates}
\|P\|_{\Omega}=\max_{x\in \Omega}\sum_{j=1}^d |\lambda_j (x)|
=\max\left \{ \sum_{j=1}^{d}
 |\beta_j| : \ \sum_{j=1}^{d} \beta_j
=1, \ y=\sum_{j=1}^{d}\beta_jy^{(j)} \in T(\Omega) \right \}.
\end{equation}
The second equality in  \eqref{norm_P_by_bary_coordinates} expresses the norm of $P$
via barycentric coordinates of the points of $T(\Omega)$ relatively  to the nondegenerate 
$(d-1)$-dimensional simplex with vertices $y^{(j)}.$

By
%$\theta_n(\Pi;\Omega)$ 
 $\theta(\Pi;\Omega)$, we  denote the minimal  norm of a projector
$P:C(\Omega)\to\Pi$:
% under the condition that the corresponding nodes belong to $\Omega:$
%$$\theta_n(\Pi;\Omega):= \min\limits_{x^{(j)}\in \Omega} \|P\|_{\Omega}.$$
$$\theta(\Pi;\Omega):= \min\limits_{x^{(j)}\in \Omega} \|P\|_{\Omega}.$$
A projector whose  norm equals  
$\theta(\Pi;\Omega)$ is called
{\it minimal}.
Also we consider the following numerical characteristic of $\Omega$ representing {\it the minimal
absorption coefficient} of this set by nondegenerate simplices having vertices in $\Omega$:
$$\xi_n(\Omega):=\min \left\{ \xi(\Omega;S): \,
S \mbox{ --- $n$-dimensional simplex,} \,
\ver(S)\subset \Omega, \, \vo(S)\ne 0\right\}.$$

Suppose $x^{(1)},\ldots,x^{(d)}\in\Omega$ is an admissible set of nodes to  interpolate
 functions from $C(\Omega)$ by  polynomials from $\Pi$,  
$P:C(\Omega)\to \Pi$ is the corresponding interpolation projector.  
Then  points $y^{(j)}=T\left(x^{(j)}\right)$ form an admissible set of nodes to interpolate
functions from  $C(T(\Omega))$ by polynomials from 
$\Pi_1\left({\mathbb R}^{d-1}\right).$
Consider the interpolation projector 
$\overline{P}:C(T(\Omega))\to \Pi_1\left({\mathbb R}^{d-1}\right)$
with the nodes  
$y^{(1)},\ldots,y^{(d)}.$ If $f\in C(\Omega),$ $g\in C(T(\Omega))$ and
$g(y)=f(x)$ while $y=T(x),$ then equalities 
$Pf\left(x^{(j)}\right)=f_j$ lead to equalities
$\overline{P}g\left(y^{(j)}\right)=g_j:=g\left(y^{(j)}\right).$ Therefore,
 interpolational polynomials $p\in\Pi$ and $q\in \Pi_1\left({\mathbb R}^{d-1}\right)$
are also connected by the relation $p(x)=q(y).$ 
Let $\lambda_j\in \Pi$ be the basic Lagrange polynomials corresponding to the projector~$P,$ $\mu_j$ be  
the basic Lagrange polynomials corresponding to the projector~$\overline{P}$ (or to the simplex $S=\conv\left(y^{(1)},\ldots,y^{(d)}\right)$ $\subset$ ${\mathbb R}^{d-1}$), then~$\mu(y)=\mu(T(x))=\lambda_j(x)$.
 Hence, 
$$\|\overline{P}\|_{T(\Omega)}=  \|P\|_\Omega. $$
By this, when estimating the norm $\|P\|_{\Omega}$, it turns out to be possible to apply geometric inequalities for the norm
of projector $\overline{P}$ under linear interpolation \linebreak  on the $(d-1)$ dimensional set $T(\Omega)$. This approach
was proposed by M.\,V. Nevskii; his results are contained in papers \cite{nevskii_mais_2008}, \cite{nevskii_mais_2009} and
monograph \cite{nevskii_monograph}. Let us present the~statements we need.

 For a projector $P:C(\Omega)\to\Pi$
with nodes $x^{(j)}$, 
\begin{equation}\label{norm_P_xi_Omega_S_ineqs}
\frac{1}{2}\left(1+\frac{1}{d-1}\right)\left(\|P\|_{\Omega}-1\right)+1
\leq \xi(T(\Omega);S)\leq  \frac{d}{2}\left(\|P\|_{\Omega}-1\right)+1, 
\end{equation}
where $S$ is the $(d-1)$-dimensional simplex with  vertices $T\left(x^{(j)}\right).$ 
Also we have
\begin{equation}\label{xi_d_minus_1_T_Omega_S_theta_n_ineqs}
\frac{1}{2}\left(1+\frac{1}{d-1}\right)\left(\theta(\Pi;\Omega)-1\right)+1
\leq \xi_{d-1}(T(\Omega))\leq
\frac{d}{2}\left(\theta(\Pi;\Omega)-1\right)+1. \end{equation}
If $\|P\|_{\Omega}\ne 1$, then \eqref{norm_P_xi_Omega_S_ineqs} can be written in the form
\begin{equation}\label{norm_P_xi_Omega_S_ineqs_modified}
\frac{1}{2}\left(1+\frac{1}{d-1}\right)
\leq \frac{\xi(T(\Omega);S)-1}{\|P\|_{\Omega}-1}\leq  \frac{d}{2}.
\end{equation}
For $\theta(\Pi;\Omega)\ne1$,
%$\theta_n(\Pi;\Omega)\ne1$ 
\eqref{xi_d_minus_1_T_Omega_S_theta_n_ineqs} is equivalent to
\begin{equation}\label{xi_d_minus_1_T_Omega_S_theta_n_ineqs_modified}
\frac{1}{2}\left(1+\frac{1}{d-1}\right)
\leq \frac{\xi_{d-1}(T(\Omega))-1}{\theta(\Pi;\Omega)-1}\leq
\frac{d}{2}. \end{equation}

We call a point $y=T(x)\in T(\Omega)$ {\it a $1$-point with respect to the simplex 
$S=\conv\left(y^{(1)},\ldots,y^{(d)}\right)$},
if
$$ \|P\|_\Omega=\sum_{j=1}^d |\lambda_j(x)|$$
and among the numbers $\lambda_j(x)$ there is the only one negative. If such a point exists, then
the right-hand inequality in \eqref{norm_P_xi_Omega_S_ineqs} becomes an equality.
 This property was proved in \cite{nevskii_mais_2008} in the equivalent form.
The notions of a $1$-vertex of a cube and a $1$-point \linebreak of an arbitrary set were introduced 
in \cite{nevskii_mais_2009} and \cite{nevskii_ukhalov_2019} correspondingly.

If a  $1$-point of 
$T(\Omega)$  exists for a simplex $S=\conv\left(T\left(x^{(1)}\right),\ldots,T\left(x^{(d)}\right)\right)$ and
$S$ satisfies
$\xi(T(\Omega);S)=\xi_{d-1}(T(\Omega))$, then the right-hand relation in   \eqref{xi_d_minus_1_T_Omega_S_theta_n_ineqs} turns to equality and the projector
$P:C(\Omega)\to \Pi$  with these nodes  $x^{(j))}$ is minimal.

The variant $\Pi=\Pi_1\left({\mathbb R}^n\right)$ was studied in the authors' works
(see, e.\,g., \cite{nevskii_mais_2009}, \cite{nevskii_monograph}, 
\cite{nevskii_ukhalov_2018},~\cite{nevskii_ukhalov_2019}).
 %In this case
 This time $d=$
$\dim \Pi_1\left({\mathbb R}^n\right)$ $=$ $n+1$ and 
the~mapping $T$ is the identity operator. The most detailed
%interesting 
results concerns the cases when $\Omega$ is an $n$-dimensional cube or an $n$-dimensional Euclidean ball.

In the present paper we discuss relations 
\eqref{norm_P_xi_Omega_S_ineqs}--\eqref{xi_d_minus_1_T_Omega_S_theta_n_ineqs_modified} for
polynomial interpolation \linebreak on a segment. Further 
 $n=1$, $\Pi=\Pi_{k}\left({\mathbb R}^1\right)$, $k\geq 1$, $d=
\dim \Pi_{k}\left({\mathbb R}^1\right)=k+1$,
$\varphi_j(x)=x^j$ $(1\leq j\leq {k})$. Let us take $\Omega=[-1,1]$, then
$$T(\Omega)=T([-1,1])=\{(x,\ldots,x^{k})\in {\mathbb R}^{k}: \, -1\leq x\leq 1\}.$$

In Sections 2--4, we consider the cases $k=2, 3, 4$. In Section 5, we present numerical  estimates
of  the  values $\theta\left(\Pi_k\left({\mathbb R}^1\right);  [-1,1]\right)$ and
   $\xi_k\left(T([-1,1])\right)$, $1\leq k\leq10$.
Finally, Sections 6--7 contain the material related to  regular and  Chebyshev nodes.
We remark all the cases when there are right-hand  equalities in  \eqref{norm_P_xi_Omega_S_ineqs}--\eqref{xi_d_minus_1_T_Omega_S_theta_n_ineqs_modified}.  

The numerical results were obtained by A.\,Yu. Ukhalov; the detailed data are~given in
%the base Mendeley Data 
\cite{ukhalov_mendeley_2023}.
% базе Mendeley Data (Ukhalov, Alexey (2023), <<Supplementary materials for the article \linebreak “On a geometric approach to the estimation of interpolation projectors"\,>>, Mendeley Data, V1, \linebreak doi: 10.17632/snh5m99yxr.1).
The computations were performed using the system %the
 Wolfram Ma\-the\-matica 
 %system of computer mathematics 
 (see,
e.\,g.,  \cite{wellin}, \cite{mangano}, \cite{wolfram}).
Specially written programs in
% the ~
C++ 
%language 
were also used.
To~invert matrices and solve extremal problems,  functions of
% the 
DLIB library \cite{king} were used.

\section{Quadratic Interpolation on a Segment}\label{sec2}
The simplest case of  nonlinear  interpolation is quadratic interpolation on a segment.
An analytical solution of the minimal projector problem with the indicated geometric approach is given in \cite{nevskii_mais_2008}, computer methods were used in
\cite{bogomolova}. Let us discuss this case as an illustration.
 
 It is well known (see, e.\,g., \cite{pashkovskij}) that  the minimal norm of an interpolation projector is realized for regular nodes and equals $\dfrac{5}{4}$. We claim that this result can be also obtained with the use of   
\eqref{norm_P_xi_Omega_S_ineqs}--\eqref{xi_d_minus_1_T_Omega_S_theta_n_ineqs}.
Additionally, it turns out that there are infinitely many minimal projectors.

\begin{figure}[h]
  \centering \includegraphics[width=14cm]
  %\textwidth]
  {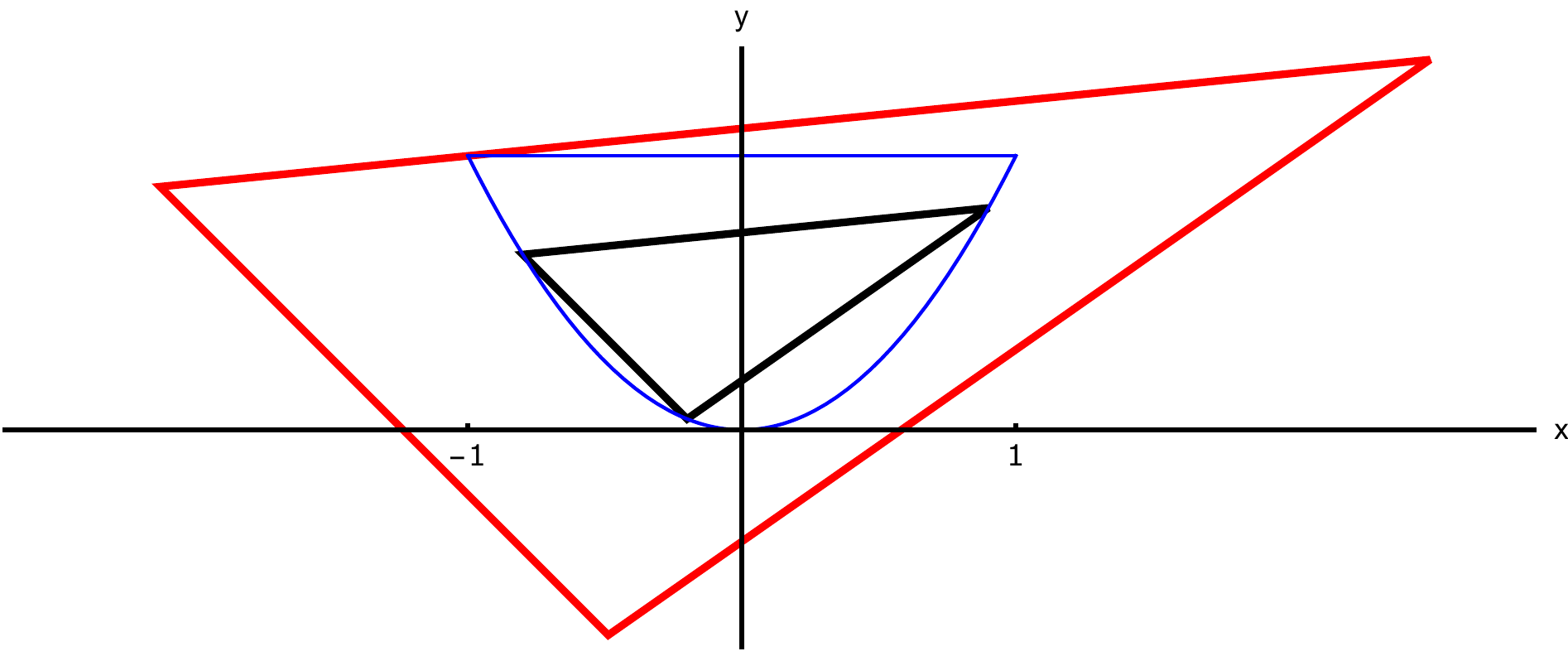}
  \caption
  {The absorption of the parabolic sector by a triangle}
  \label{fig:Triangle}
\end{figure}

If $
%\Omega=[-1,1],
 \Pi=\Pi_2\left({\mathbb R}^1\right)$, then $d=\dim \Pi=3,$ i.\,e., $k=d-1=2.$  
The mapping $T$ has the form $x\longmapsto (x,x^2)$ and the set 
$T(\Omega)=T([-1,1])=\{(x,x^2)\in {\mathbb R}^2: \, -1\leq x\leq 1\}$ \linebreak
is a piece of parabola. For the interpolation  nodes
%$-1\leq x^{(1)}<x^{(2)}<x^{(3)}\leq 1$ 
$-1\leq r<s<t\leq 1$, we have
$${\bf A} =
\left( \begin{array}{ccc}
1&r&r^2\\
1&s&s^2\\
1&t&t^2
\end{array}
\right),
$$  
the simplex $S$ is a triangle with  the vertices  
$(r,r^2), \,(s,s^2), \,(t,t^2)$ lying on this part of parabola.
The absorption of the parabolic sector by this triangle is shown in~Figure~\ref{fig:Triangle}\,.

The convexity of the function $\psi(x)=x^2$ implies that a $1$-point of the set  $T([-1,1])$ 
with respect to $S$ exists for any nodes.
Hence, there is an equality in~the~right-hand part of \eqref{norm_P_xi_Omega_S_ineqs}:
\begin{equation}\label{norm_P_xi_eq_for_quadratic_interp}
\xi(T(\Omega);S)= \frac{3}{2}\left(\|P\|_{\Omega}-1\right)+1=\frac{3\|P\|_{\Omega}-1}{2}. 
\end{equation}
Since  \eqref{norm_P_xi_eq_for_quadratic_interp} holds true for an arbitrary projector 
$P:C[-1,1]\to 
\Pi_2({\mathbb R}^1)$, we have  the~right-hand equality also in  
\eqref{xi_d_minus_1_T_Omega_S_theta_n_ineqs}:
$$\xi_2(T(\Omega))
=\frac{3\theta\left(\Pi; \Omega\right)-1}{2}.$$

So, in the quadratic case finding the minimum norm of a projector
$\theta\left(\Pi; \Omega\right)$
is~equivalent to calculating 
$\xi_2(T(\Omega))$,
 i.\,e., the minimum absorption coefficient  of the specified part of the parabola by a triangle. 
 As nodes of a minimum projector, one must take the first coordinates of vertices of the result triangle.
 Technically, the~problem can be reduced to a~triangle $S$ with vertices $(-r,r^2),$ $(0,0),$ $(r,r^2),$
$0<r\leq 1.$
For~this triangle, the extremal points 
%$y\in T([-1,1])$ 
$y\in T(\Omega)$ can be only
$\displaystyle(\pm 1,1), \left(\pm \frac{r}{2}, \frac{r^2}{4}\right).$ At~each of the latter two points
a tangent line to the parabola is parallel to a side of~$S.$ Calculations give
$$\xi(T(\Omega);S)=\max\left(\frac{11}{8},\frac{3}{r^2}-2\right), \quad
\|P\|_{\Omega}=\max\left(\frac{5}{4},\frac{2}{r^2}-1\right).$$
If $$\frac{2\sqrt{2}}{3}=0.942809\ldots\leq r\leq 1,$$ then 
these values do not depend on 
$r$ and are equal correspondingly to  \,$\dfrac{11}{8}$ and \,$\dfrac{5}{4}$\,;
these are minimal possible values.
Thus, 
$$\xi_2(T(\Omega))=\frac{11}{8}, \quad \theta\left(\Pi; \Omega\right)=\frac{5}{4}.$$
The minimal is any projector with vertices $-r,$ $0,$ $r$ provided $r\in\left[\dfrac{2\sqrt{2}}{3},1\right]$
and there are no other minimal projectors.
Note that the relations
\eqref{xi_d_minus_1_T_Omega_S_theta_n_ineqs}--\eqref{xi_d_minus_1_T_Omega_S_theta_n_ineqs_modified} 
look like
$$\frac{11}{16}<\frac{11}{8}=\frac{11}{8}, \quad \frac{3}{4}<\frac{3}{2}=\frac{3}{2}.$$

%%%%%%%%%%%%%%%%%%%%%%%%%%%%%%%%%%%%%%%%%%%%%%%
\section{Cubic Interpolation on a Segment}
In the case $\Pi=\Pi_3\left({\mathbb R}^1\right)$  we have $d=\dim \Pi=4,$ $k=d-1=3.$  
Since $T$ has the form $x\longmapsto (x,x^2,x^3),$ then 
$$T(\Omega)=T([-1,1])=\{(x,x^2,x^3)\in {\mathbb R}^3: \, -1\leq x\leq 1\}.$$
This is a three-dimensional line with  endpoints 
 $(-1,1,-1)$ and $(1,1,1)$
  whose projections on the coordinate planes are
 % which has projections upon the coordinate planes
 congruent to the curves $Y=X^2,$ $Y=X^3$, and 
$X(t)=t^2,$ $Y(t)=t^3;$ the latter one has the zero angle at the point
$X=Y=0.$

\begin{figure}[h]
  \centering \includegraphics[width=12cm]{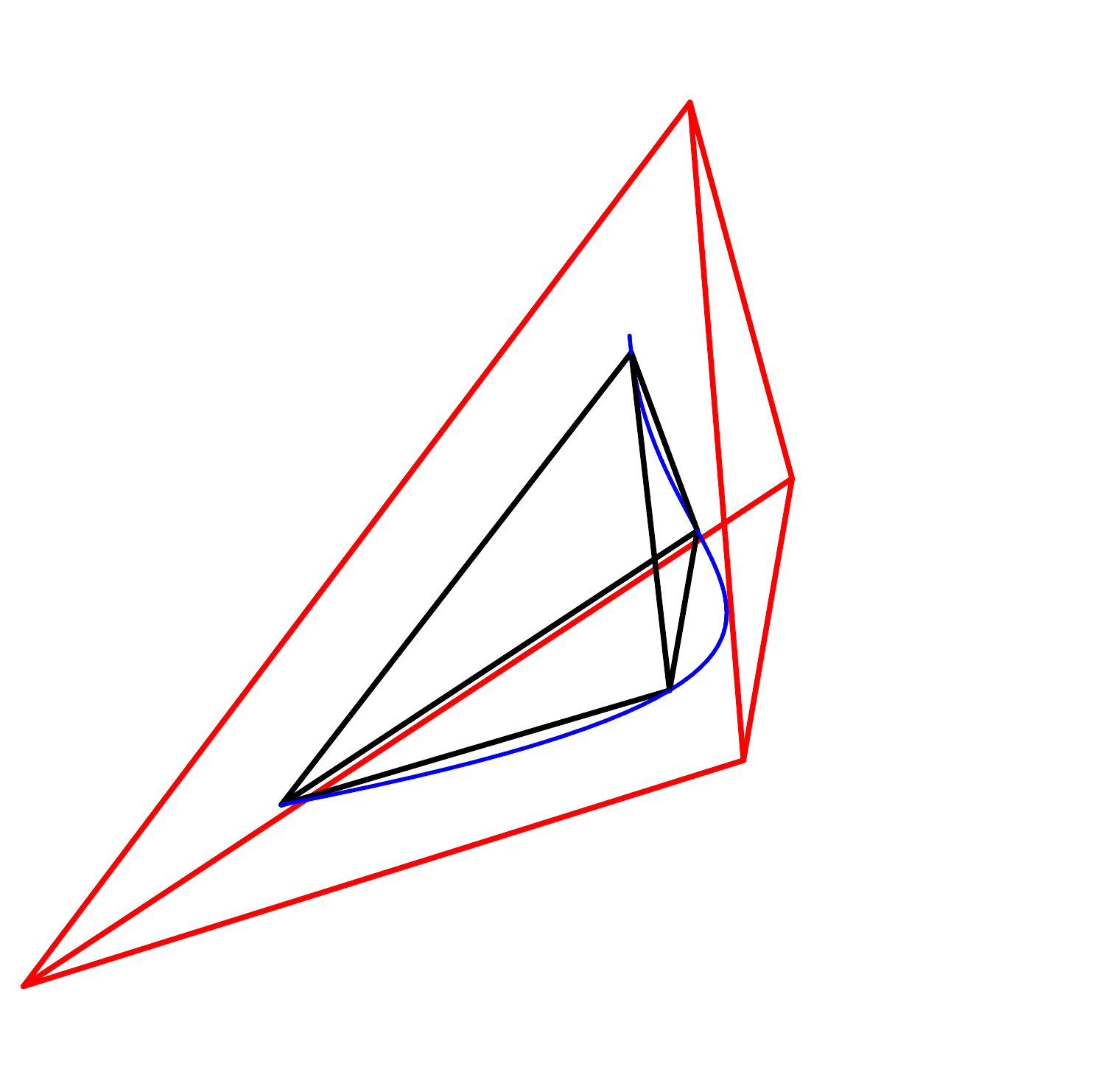}
\caption{The absorption of the set $T([-1,1])$  by a tetrahedron  }
 \label{fig:Zagogulina}
\end{figure}

 For the nodes 
$-1\leq r<s<t<u\leq 1$, 
$${\bf A} =
\left( \begin{array}{cccc}
1&r&r^2&r^3\\
1&s&s^2&s^3\\
1&t&t^2&t^3\\
1&u&u^2&u^3\\
\end{array}
\right),
$$  
the simplex $S$ is a tetrahedron with the vertices 
\begin{equation}\label{nodes_cubic_interpolation}
(r,r^2,r^3), \quad (s,s^2,s^3), \quad (t,t^2,t^3), \quad (u, u^2, u^3)
\end{equation}
 containing in $T(\Omega)$. The absorption of the set  $T(\Omega)$ by this tetrahedron  is shown 
 in~ Figure~\ref{fig:Zagogulina}\,.
  
Minimal values of $\|P\|_\Omega$ and $\xi(T(\Omega);S)$  were found using a computer. 
The minimum of $\|P\|_\Omega$ equal to $1.422919...$  is delivered by symmetrical nodes
\begin{equation}\label{min_projector_nodes_cubic}
 -1 \,, \quad -0.417791... \, , \quad   0.417791... \, , \quad  1 \,.
 \end{equation}
The minimum of  $\xi(T(\Omega);S)$ equal to $1.635778...$ is delivered by a tetrahedron
whose vertices appear after applying the mapping $T$ to the points 
$$ -1, \quad -0.481618... \, , \quad  0.481618... \, , \quad \  1 \,.$$

 So,  the computer calculations give
$\theta\,(\Pi;\Omega)=1.422919... \,,$  \, $\ \xi_3(T(\Omega))=1.635778... \,.$
With such values both the inequalities in \eqref{xi_d_minus_1_T_Omega_S_theta_n_ineqs} are strict and have the
form $1.28194... < 1.63577... < 1.84583... \,.$ Relation \eqref{xi_d_minus_1_T_Omega_S_theta_n_ineqs_modified} may be written as
$\frac{2}{3}<1.503307\ldots  < 2.$

Under quadratic interpolation, a $1$-point of  $T(\Omega)$  
exists for an arbitrary projector.  But  in the case considered this is not so.   For example, let us take 
$$r=-\frac{\sqrt{2+\sqrt{2}}}{2}, \quad s
=-\frac{\sqrt{2-\sqrt{2}}}{2}, \quad
t=\frac{\sqrt{2-\sqrt{2}}}{2}, \quad
u=\frac{\sqrt{2+\sqrt{2}}}{2}.$$
These are the zeros of the Chebyshev polynomial of degree 4, i.\,e.,  
%$T_4(x)=
$8x^4-8x^2+1$; we~call them  {\it Chebyshev nodes.} Then
$\|P\|_{\Omega}= \sqrt{2+\sqrt{2}} = 1.847759\ldots \, , \
 \xi(T(\Omega); S)=2.496605\ldots \,  
$,
and inequalities \eqref{norm_P_xi_Omega_S_ineqs} have the form $1.5651727\ldots<2.496605\ldots$ \linebreak
$<2.695518\ldots \ .$
Since the right-hand relation in \eqref{norm_P_xi_Omega_S_ineqs} is not an equality,
a  $1$-point of 
$T(\Omega)$ with respect to the simplex with vertices  \eqref{nodes_cubic_interpolation} 
does not exist.

Nethertheless, our calculations show  that for some nodes the right-hand relation in \eqref{norm_P_xi_Omega_S_ineqs}  turns an equality  (with  the accurancy $10^{-11}$). This takes place
for  regular nodes and also for the nodes \eqref{min_projector_nodes_cubic} of the minimal projector.

For regular nodes,
%$-1,-\dfrac{1}{3},\dfrac{1}{3},1$ 
$\|P\|_{\Omega} = 1.63113030\ldots \, , \
 \xi(T(\Omega); S)=2.26226061\ldots \,$; relation   \eqref{norm_P_xi_Omega_S_ineqs} takes the form 
 $1.420753\ldots < 2.262260\ldots = 2.262260\ldots$\,. 
 Remark that we maid the calculations of $\|P\|_{\Omega}$ и $\xi(T(\Omega); S)$ with the use of two various programs. The~latter equality holds true with the accurancy not less than   $10^{-18}$.

 For the nodes \eqref{min_projector_nodes_cubic},
the projector norm is minimal: $\|P\|_{\Omega}=\theta\,(\Pi;\Omega)=1.422919\ldots,$ the absorption
coefficient is
 $\xi(T(\Omega); S)=1.845839\ldots$, and relation \eqref{norm_P_xi_Omega_S_ineqs}
 is written as~$1.281946\ldots < 1.845839\ldots = 1.845839\ldots$ .

Of course, approximate calculations, even  while having good accurancy, do not garantee the existence
of the right-hand equality in  \eqref{norm_P_xi_Omega_S_ineqs}.
The strong proof gives an~approach with the use of the notion of a $1$-point (see Section \ref{nev_s1}). 

For  regular nodes
$$x^{(1)}=-1,\quad x^{(2)}=-\dfrac{1}{3}, \quad x^{(3)}=\dfrac{1}{3}, \quad x^{(4)}=1$$
and for the point $x^*= -0.699055\ldots$,  we have
$\|P\|_\Omega= 1.631130\ldots=\sum |\lambda_j\left(x^*\right)|$, while
$$\lambda_1\left(x^*\right)= 0.890801\ldots,  \quad  \lambda_2\left(x^*\right)= -0.315565\ldots, $$ $$ \lambda_3\left(x^*\right)=0.360848\ldots,
\quad  \lambda_4\left(x^*\right)= 0.063915\ldots .$$
Consequently, $y^*=T\left(x^*\right)$ is a  $1$-point of  $T(\Omega)$ with respect to a
tetrahedron with  vertices
$T\left(x^{(j)}\right)$. 
In accordance with  Section \ref{nev_s1}, the right-hand relation
 in~\eqref{norm_P_xi_Omega_S_ineqs} turns into equality. Note that $\det({\bf A})= 1.053497\ldots .$

For  nodes 
$$x^{(1)}=-1,\quad x^{(2)}= -0.417791\dots, \quad x^{(3)}= 0.417791\dots, \quad x^{(4)}=1$$
and for $x^{**}=  -0.733172\ldots$, we have 
$\|P\|_\Omega= 1.422919\ldots=\sum |\lambda_j\left(x^{**}\right)|$.
Since 
$$\lambda_1\left(x^{**}\right)=  -0.211459\ldots,  \quad  \lambda_2\left(x^{**}\right)=  0.771708\ldots, $$ $$ \lambda_3\left(x^{**}\right)=0.381082\ldots,
\quad  \lambda_4\left(x^{**}\right)= 0.058668\ldots ,$$
 $y^{**}=T\left(x^{**}\right)$ is a $1$-point of  $T(\Omega)$ with respect to a tetrahedron with
vertices 
$T\left(x^{(j)}\right)$. The existence of a 1-point means that also in this case we have
an equality in the right-hand part of 
\eqref{norm_P_xi_Omega_S_ineqs}.   Here $\det({\bf A})=  -1.138679\ldots .$

 \section{Interpolation by Polynomials from $\Pi_4\left({\mathbb R}^1\right)$}\label{sec4}
 
  If
 $\Pi=\Pi_4\left({\mathbb R}^1\right)$,  then  $d=\dim \Pi=5,$ $k=d-1=4.$  
The mapping $T$ takes the~form $x\longmapsto (x,x^2,x^3,x^4),$   therefore

$$T(\Omega)=T([-1,1])=\{(x,x^2,x^3,x^4)\in {\mathbb R}^4: \, -1\leq x\leq 1\}.$$
 For nodes
$-1\leq r<s<t<u<v\leq 1$, we have 
$${\bf A} =
\left( \begin{array}{ccccc}
1&r&r^2&r^3&r^4\\
1&s&s^2&s^3&s^4\\
1&t&t^2&t^3&t^4\\
1&u&u^2&u^3&u^4\\
1&v&v^2&v^3&v^4\\
\end{array}
\right).
$$  
Coordinates of vertices of  simplex  $S\subset {\mathbb R}^4$ are written in the rows of matrix ${\bf A}$
starting from the second column.

 The Chebyshev nodes, i.\,e., the zeroes of the $5$th degree Chebyshev
 polynomial  
%$T_5(x)=
$16 x^5-20 x^3+5 x$, are
$$r=-\frac{\sqrt{5+\sqrt{5}}}{2\sqrt{2}}, \quad s=-\frac{\sqrt{5-\sqrt{5}}}{2\sqrt{2}}, \quad t=0, \quad 
u=\frac{\sqrt{5-\sqrt{5}}}{2\sqrt{2}}, \quad v=\frac{\sqrt{5+\sqrt{5}}}{2\sqrt{2}}.$$
In this case,
$\|P\|_\Omega=\dfrac{1+4\sqrt{5}}{5}=1.988854\ldots \, ,$ but this is not the possible minimum value.

Let us give the results of numerical calculations.
The minimal $\|P\|_\Omega$ equal to $1.559490...$ is achieved at symmetrical nodes
$$ \-1 \, , \quad  -0.620911... \, ,  \quad 0 \, ,  \quad  0.620911... \, , \quad 1 \, .$$
The minimal $\xi(T(\Omega);S)$ equal to $1.981193...$ is delivered by a simplex with vertices constructed
by operator $T$ from the points
$$ -1 \, , \quad -0.650738... \, ,\quad   0 \, ,\quad  0.650738... \, , \quad
1 \, .$$
Thus,
$\theta\left(\Pi;\Omega\right)=1.559490... \,, \ \xi_4(T(\Omega))=1.981193...\,.$
Inequalities in \eqref{xi_d_minus_1_T_Omega_S_theta_n_ineqs} 
are both strict and have the form $
1.349681... < 1.981193... < 2.398725... \,.$ 

%%%%%%%%%%%%%%%%%%%%%%%%%%%%%%%%%

   \section{Estimates of $\theta\left(\Pi_k\left({\mathbb R}^1\right);  [-1,1]\right)$  and 
   $\xi_k\left(T([-1,1])\right)$ for $1\leq k\leq10$}\label{sec5}

  In this section, we present numerical estimates for
   minimal norms of interpolation projections acting from $C[-1,1]$ on 
   $\Pi_k\left({\mathbb R}^1\right)$ 
   %for $1\leq k\leq10$
   and minimum absorption coefficients for~$1\leq k\leq10$.
 %  in the same cases.
   For brevity, we write
    $\theta_k:=\theta\left(\Pi_k\left({\mathbb R}^1\right); [-1,1]\right)$, \, $\xi_k:=\xi_k\left(T ([-1,1])\right)$.

  \begin{table}[h!]
	  \caption{Minimal norms of projectors  and~corresponding absorption coefficients}
	\label{table:1}
	\centering
	\medskip
  \begin{tabular}{|c|c|c|}
  	\hline
  	$k$ & $\theta_k\leq$ & $\xi(S)$ \\
  	\hline
	1 & 1& 1\\
	\hline
2 & 1.25& 1.375\\
\hline
3 & 1.422919\dots & 1.845839\dots \\
\hline
4 & 1.559490\dots & 2.224196\dots \\
\hline
5 & 1.672210\dots & 2.574785\dots \\
\hline
6 & 1.768134\dots & 2.911143\dots \\
\hline
7 & 1.851599\dots & 3.239031\dots \\
\hline
8 & 1.925457\dots & 3.561425\dots \\
\hline
9 & 1.991685\dots & 3.880036\dots \\
\hline
10 & 2.051705\dots & 4.195926\dots \\
   	\hline
   \end{tabular}
\end{table}
\medskip

    \begin{table}[h!]
    	\caption{Minimal absorption coefficients  and~corresponding norms of projectors}
    	\label{table:2}
    	\centering
    	\medskip
    	\begin{tabular}{|c|c|c|}
    		  	\hline
$k$ & $\xi_k\leq$ & $||P||$ \\
   	\hline    		
 	1 & 1& 1\\
 	\hline   
2 & 1.375000\dots & 1.250000\dots \\
\hline   
3 & 1.635778\dots & 1.604018\dots \\
\hline   
4 & 1.981193\dots & 1.626067\dots \\
\hline   
5 & 2.210535\dots & 1.782786\dots \\
\hline   
6 & 2.455130\dots & 1.858521\dots \\
\hline   
7 & 2.678509\dots & 1.962845\dots \\
\hline   
8 & 2.907301\dots & 2.029565\dots \\
\hline   
9 & 3.128316\dots & 2.108072\dots \\
\hline   
10 & 3.351866\dots & 2.164915\dots \\
   	\hline    		
    	\end{tabular}
    \end{table}
    \medskip

   Table \ref{table:1} contains the number $k$ (the degree of an interpolation polynomial), the~upper bound for~$\theta_k$, and the value $\xi(S):=\xi(T([-1,1]);S)$, i.\,e., the absorption coefficient of the set 
   $T([-1,1])$ by the simplex with vertices  $T\left(x^{(j)}\right)$, where $x^{(j)}$ are 
   the~nodes of interpolation projector delivering the estimate of $\theta_k$. The values $\theta_1=1$ and $\theta_2=1.25$
   are precise. The right-hand equality in  \eqref{norm_P_xi_Omega_S_ineqs} takes places for  
   $k=1,2,3$.
 
 In Table  \ref{table:2}, we give the upper bounds for minimal absorption coefficients
 $\xi_k$ and~the~norms
$\|P\|:=\|P\|_{[-1,1]}$ of  projectors with the nodes constructed by vertices of extremal simplices. The values $\xi_1=1$ and $\xi_2=1.375$ are precise.  The~second equality in  \eqref{norm_P_xi_Omega_S_ineqs} takes places for  
   $k=1,2$.

The  numerically found nodes of minimal projections  and the parametric coordinates of the vertices of extremal simplices are given in \cite{ukhalov_mendeley_2023}. There one can also find the estimates for $\theta_k$ and $\xi_k$ with a~large number of decimal signs. \linebreak  See https://doi.org/10.17632/snh5m99yxr.1.

%%%%%%%%%%%%%%%%%%%%%%%%%%%%  
   
  \section{Regular Nodes}\label{sec6}

  Computer calculations of absorption coefficients
    and projector norms for regular nodes (see Table \ref{table:3}) are complicated by a rapid decrease in the modulus
   of  determinant of matrix ${\bf A}$. For example, for $k=11$ we have $\det ({\bf A})= 3.63581\cdot 10^{-14}$. Further, as $k$ increases, the absolute value of the determinant rapidly decreases.
    When computing on systems with standard floating point representation, it is difficult to guarantee the accuracy of the results.
    For this reason, we present the values of the considered quantities only for $1\leq k \leq 10$. Note that the right-hand equality in \eqref{norm_P_xi_Omega_S_ineqs} holds for
    $k=1,2,3$.

       \begin{table}[h!]
   	\caption{Absorption coefficients and norms of~projectors for regular nodes}
   	\label{table:3}
   	\centering
   	\medskip
   	\begin{tabular}{|c|c|c|}
   		\hline
   		$k$ & $\xi(S)$ & $||P||$ \\
   		   		\hline
   	1 & 1 & 1\\
   	   		\hline
2 & 1.375& 1.25\\
   		\hline
3 & 2.262260\dots & 1.631130\dots \\
   		\hline
4 & 3.812500\dots & 2.207824\dots \\
   		\hline
5 & 6.167317\dots & 3.106301\dots \\
   		\hline
6 & 9.461457\dots & 4.549341\dots \\
   		\hline
7 & 13.824447\dots & 6.929739\dots \\
   		\hline
8 & 21.876588\dots & 10.945645\dots \\
   		\hline
9 & 41.283675\dots & 17.848612\dots \\
   		\hline
10 & 72.576233\dots & 29.899955\dots \\
   		\hline    		
   	\end{tabular}
   \end{table}
  \medskip

 \begin{table}[h!]
   	\caption[Коэффициенты поглощения  и~нормы проекторов для чебышёвских узлов]{Absorption coefficients and norms of~projectors for Chebyshev nodes}
   	\label{table:4}
   	\centering
   	\medskip
   	\begin{tabular}{|c|c|c|}
   		\hline
   		$k$ & $\xi(S)$ & $||P||$ \\
   		\hline
1 & 1.414213\dots & 1.414213\dots \\
\hline
2 & 2.000000\dots & 1.666666\dots \\
\hline
3 & 2.496605\dots & 1.847759\dots \\
\hline
4 & 2.962610\dots & 1.988854\dots \\
\hline
5 & 3.414213\dots & 2.104397\dots \\
\hline
6 & 3.857835\dots & 2.202214\dots \\
\hline
7 & 4.296558\dots & 2.287016\dots \\
\hline
8 & 4.732050\dots & 2.361856\dots \\
\hline
9 & 5.165299\dots & 2.428829\dots \\
\hline
10 & 5.596925\dots & 2.489430\dots \\
\hline
11 & 6.027339\dots & 2.544766\dots \\
\hline
12 & 6.456823\dots & 2.595678\dots \\
   		\hline    		
   	\end{tabular}
   \end{table}
   \medskip
   
   \section{Chebyshev Nodes }\label{sec7}
   
The table \ref{table:4} shows the results of calculating absorption coefficients and projector norms for Chebyshev nodes
 (i.\,e., the nodes coinciding with the zeros of Chebyshev polynomial of the required degree). The right-hand equality in \eqref{norm_P_xi_Omega_S_ineqs} holds when \linebreak $k=1,2$.

As in the case of regular nodes, one can vouch for the accuracy of calculations of~absorption coefficients  only for small $k$. We give values for $1\leq k \leq 12$.
Already for $k=12$ we have $\det ({\bf A}) = 3.68529\cdot 10^{-15}$. As $k$ increases, $|\det ({\bf A})|$   becomes even smaller, which can lead to a decrease in precision when using C++ variables of~the~double type.
Projector norms for Chebyshev nodes were calculated  using the~exact formula
(see, e.\,g.,~\cite{pashkovskij}), so this remark does not apply to the values of~the~norms.

Note that the C++ program we use to calculate the norms of projections gives a satisfactory agreement 
%between the values and 
with the exact formula at least for $k\leq 30$.

%\clearpage

\end{document}